\documentclass[11pt]{article}

\usepackage{amssymb}
\usepackage[english]{babel}
\usepackage{graphicx}

\newtheorem{thm}{Theorem}[section]
\newtheorem{pro}[thm]{Proposition}
\newtheorem{defi}[thm]{Definition}

\def\G{\Gamma}
\def\dist{\mathop{\rm dist }\nolimits}

\begin{document}

\title{Dual Concepts of Almost Distance-Regularity
and the Spectral Excess Theorem\thanks{This version is published in Discrete
Math. 312 (2012), 2730--2734. Research supported by the Ministerio de
Educaci\'on y Ciencia (Spain) and the European Regional Development Fund under
project MTM2008-06620-C03-01, and by the Catalan Research Council under project
2009SGR1387.}}

\author{C. Dalf\'{o}$^\dag$, E.R. van Dam$^\ddag$, M.A. Fiol$^\dag$, E.
Garriga$^\dag$
\\ \\
{\small $^\dag$Universitat Polit\`ecnica de Catalunya, Dept. de Matem\`atica Aplicada IV} \\
{\small Barcelona, Catalonia} {\small (e-mails: {\tt
\{cdalfo,fiol,egarriga\}@ma4.upc.edu})} \\
{\small $^\ddag$Tilburg University, Dept. Econometrics and O.R.} \\
{\small Tilburg, The Netherlands} {\small (e-mail: {\tt
edwin.vandam@uvt.nl})} \\
}

\date{}

\maketitle

\begin{abstract}
Generally speaking, `almost distance-regular' graphs share some, but not
necessarily all, of the regularity properties that characterize
distance-regular graphs. In this paper we propose two new dual concepts of
almost distance-regularity, thus giving a better understanding of the
properties of distance-regular graphs. More precisely, we characterize
$m$-partially distance-regular graphs and $j$-punctually eigenspace
distance-regular graphs by using their spectra. Our results can also be seen as
a generalization of the so-called spectral excess theorem for distance-regular
graphs, and they lead to a dual version of it.
\end{abstract}



\section{Preliminaries}
Almost distance-regular graphs, recently studied in the literature, are graphs
which share some, but not necessarily all, of the regularity properties that
characterize distance-regular graphs. Two examples of the former are partially
distance-regular graphs \cite{p91} and $m$-walk-regular graphs \cite{dfg09}.

In this paper we propose and characterize two dual concepts of almost
distance-regularity, and study some cases where distance-regularity is
attained. As in the theory of distance-regular graphs, the two proposed
concepts lead to several duality results. Our results can also be seen as a
generalization of the so-called spectral excess theorem for distance-regular
graphs (see \cite{fg97}; for short proofs, see \cite{vd08,fgg09}). This theorem
characterizes distance-regular graphs by their spectra and the average number
of vertices at extremal distance. A dual version of this theorem is also
derived.

We use standard concepts and results for distance-regular graphs
\cite{biggs,bcn}, spectral graph theory \cite{cds80,g93}, and spectral and
algebraic characterizations of distance-regular graphs \cite{f02}. Moreover,
for some more details and other concepts of almost distance-regularity  (such
as distance-polynomial and partially distance-regular graphs), we refer the
reader to our recent paper \cite{ddfgg10}. In what follows, we recall the main
concepts, terminology, and results involved.

Let $\Gamma$ be a simple, connected, $\delta$-regular graph, with vertex set
$V$, order $n=|V|$, and adjacency matrix $\textbf{\emph{A}}$. The {\it
distance} between two vertices $u$ and $v$ is denoted by $\dist (u,v)$, so the
{\it diameter} of $\Gamma$ is $D=\textrm{max}_{u,v\in V}\dist(u,v)$. The set of
vertices at distance $i$ from a given vertex $u\in V$ is denoted by
$\Gamma_i(u)$, for $i=0,1,\ldots,D$. The {\em distance-$i$ graph} $\G_i$ is the
graph with vertex set $V$ and where two vertices $u$ and $v$ are adjacent if
and only if $\dist(u,v)=i$ in $\G$. Its adjacency matrix $\textbf{\emph{A}}_i$
is usually referred to as the {\em distance-$i$ matrix} of $\G$. The spectrum
of $\G$ is denoted by $ \textrm{sp}\,\G =
\{\lambda_0^{m_0},\lambda_1^{m_1},\ldots, \lambda_d^{m_d}\}, $ where the
different eigenvalues of $\G$ are in decreasing order,
$\lambda_0>\lambda_1>\cdots >\lambda_d$, and the superscripts stand for their
multiplicities $m_i=m(\lambda_i)$.

\subsection{The predistance and preidempotent polynomials}

From the spectrum of $\Gamma$, we consider the {\em
predistance polynomials} $\{p_i\}_{0\le i\le d}$ which are orthogonal with respect to the
following scalar product in $\mathbb{R}_d[x]$:
\begin{equation}
\label{product}
\langle f, g\rangle_{\vartriangle} =\frac{1}{n}\textrm{tr}\,
(f(\textbf{\emph{A}})g(\textbf{\emph{A}}))=\frac{1}{n} \sum_{i=0}^d m_i f(\lambda_i) g(\lambda_i),
\end{equation}
and which satisfy $\textrm{deg}\,p_i=i$ and $\langle p_i,p_j
\rangle_\vartriangle= \delta_{ij}p_i(\lambda_0)$, for all $i,j=0,1,\ldots,d$.
For more details, see \cite{fg97}. Like every sequence of orthogonal
polynomials, the predistance polynomials satisfy a three-term recurrence of the
form
\begin{equation}
\label{recur-pol}
xp_i=\beta_{i-1}p_{i-1}+\alpha_i p_i+\gamma_{i+1}p_{i+1},\qquad i=0,1,\ldots,d,
\end{equation}
with $\beta_{-1}=\gamma_{d+1}=0$. Some basic properties of these coefficients,
such as $\alpha_i+\beta_i+\gamma_i=\lambda_0$ for $i=0,1,\ldots, d$, and
$\beta_i n_i=\gamma_{i+1}n_{i+1}\neq0$ for $i=0,1,\ldots, d-1$, where
$n_i=\|p_i\|_{\vartriangle}^2=p_i(\lambda_0)$, can be found in \cite{cffg09}.
Let $\omega_i$ be the leading coefficient of $p_i$. Then, from the above
recurrence and since $p(0)=1$, it is immediate that $\omega_i=
(\gamma_1\gamma_2\cdots \gamma_i)^{-1}$ for $i=1,\ldots,d$.

For any graph, the sum of all the predistance polynomials gives the {\em
Hoffman polynomial} $H$ satisfying $H(\lambda_i)=n\delta_{0i}$,
$i=0,1,\ldots,d$, which characterizes regular graphs via the condition
$H(\textbf{\emph{A}})=\textbf{\emph{J}}$, the all-$1$ matrix \cite{hof63}. Note
that the leading coefficient $\omega_d$ of $H$ (and also of $p_d$) is
$\omega_d=n/\pi_0$.

From the predistance polynomials, we define the so-called  {\em preidempotent
polynomials} $q_j$, $j=0,1,\ldots, d$, by
$$
q_j(\lambda_i)= \frac{m_j}{n_i}p_i(\lambda_j),\qquad i=0,1,\ldots, d,
$$
which are orthogonal with respect to the scalar product
\begin{equation}
\label{product-preadj}
\langle f,g\rangle_\blacktriangle =\frac{1}{n}\textrm{tr}\,
(f\{\textbf{\emph{A}}\}g\{\textbf{\emph{A}}\})=\frac{1}{n}\sum_{i=0}^d  n_i f(\lambda_i) g(\lambda_i),
\end{equation}
where $f\{\textbf{\emph{A}}\}=\frac{1}{\sqrt{n}}\sum_{i=0}^d f(\lambda_i)
p_i(\textbf{\emph{A}})$. Note that, since $q_j(\lambda_0)=m_j$, the duality
between the two scalar products (\ref{product}) and (\ref{product-preadj}) and
their associated polynomials is made apparent by writing
\begin{eqnarray}
\langle p_i, p_j\rangle_\vartriangle &=& \frac{1}{n}\sum_{l=0}^d m_l p_i(\lambda_l) p_j(\lambda_l)=\delta_{ij} n_i,\qquad i,j=0,1,\ldots,d, \label{basic-predistance} \\
\langle q_i, q_j\rangle_{\blacktriangle} &=& \frac{1}{n}\sum_{l=0}^d n_l q_i(\lambda_l) q_j(\lambda_l)=\delta_{ij} m_i,\qquad i,j=0,1,\ldots,d. \label{basic-preadj}
\end{eqnarray}

\subsection{Vector spaces, algebras and bases}
\label{subsec_alg}

Let $\G$ be a graph with diameter $D$, adjacency matrix $\textbf{\emph{A}}$ and
$d+1$ distinct eigenvalues. We consider the vector spaces ${\cal A}=
\mathbb{R}_{d}[\textbf{\emph{A}}] = \linebreak \textrm{span}
\{\textbf{\emph{I}}, \textbf{\emph{A}}, \textbf{\emph{A}}^2, \ldots,
\textbf{\emph{A}}^{d}\}$ and ${\cal D}= \textrm{span}
\{\textbf{\emph{I}},\textbf{\emph{A}},\textbf{\emph{A}}_2,\ldots,\textbf{\emph{A}}_D\}$,
with dimensions $d+1$ and $D+1$, respectively. Then,  ${\cal A}$ is an algebra
with the ordinary product of matrices, known as the {\it adjacency algebra},
with orthogonal bases
$A_p=\{p_0(\textbf{\emph{A}}),p_1(\textbf{\emph{A}}),p_2(\textbf{\emph{A}}),\ldots,
p_d(\textbf{\emph{A}})\}$ and
$A_\lambda=\{\textbf{\emph{E}}_0,\textbf{\emph{E}}_1,\ldots,
\textbf{\emph{E}}_d\}$, where the matrices $\textbf{\emph{E}}_i$,
$i=0,1,\ldots,d$, corresponding to the orthogonal projections onto the
eigenspaces, are the {\it $($principal\/$)$ idempotents} of
$\textbf{\emph{A}}$. Besides, since
$\textbf{\emph{I}},\textbf{\emph{A}},\textbf{\emph{A}}^2,\ldots,\textbf{\emph{A}}^D$
are linearly independent, we have that $\textrm{dim}\, \mathcal{A}=d+1\ge D+1$
and, therefore, we always have $D\le d$ \cite{biggs}. Moreover, ${\mathcal D}$
forms an algebra with the entrywise or Hadamard product of matrices, defined by
$(\textbf{\emph{X}}\circ\textbf{\emph{Y}})_{uv}=\textbf{\emph{X}}_{uv}\textbf{\emph{Y}}_{uv}$.
We call ${\mathcal D}$ the {\em distance $\circ$-algebra}, which has orthogonal
basis $D_{\lambda}=
\{\textbf{\emph{I}},\textbf{\emph{A}},\textbf{\emph{A}}_2,\ldots,\textbf{\emph{A}}_d\}$.

From now on, we work with the vector space ${\cal T}={\cal A}+{\cal D}$, and
relate the distance-$i$ matrices $\textbf{\emph{A}}_i \in {\mathcal D}$ to the
matrices $p_i(\textbf{\emph{A}}) \in {\mathcal A}$. Note that
$\textbf{\emph{I}}$, $\textbf{\emph{A}}$, and $\textbf{\emph{J}}$ are matrices
in ${\cal A}\cap{\cal D}$ since $\textbf{\emph{J}}=H(\textbf{\emph{A}})\in
\mathcal{A}$. Recall that ${\mathcal A}={\mathcal D}$ if and only if $\G$ is
distance-regular (see \cite{biggs,bcn}). In this case, we have $D=d$, and the
predistance polynomials become the {\em distance polynomials} satisfying
$\textbf{\emph{A}}_i=p_i(\textbf{\emph{A}})$. In ${\cal T}$, we consider the
following scalar product:
\begin{equation}
\label{equationscalarproduct}
\langle\textbf{\emph{R}},\textbf{\emph{S}}\rangle=
\frac 1n\textrm{tr}\, (\textbf{\emph{RS}})= \frac
1n\textrm{sum}\,(\textbf{\emph{R}}\circ\textbf{\emph{S}}),
\end{equation}
where $\textrm{sum}\,(\textbf{\emph{M}})$ denotes the sum of all entries of
$\textbf{\emph{M}}$. Observe that the factor $1/n$ assures that
$\|\textbf{\emph{I}}\|^2=1$, whereas $\|\textbf{\emph{J}}\|^2=n$. Note also
that the {\em average degree} of $\G_i$ is
$\overline{\delta}_i=\|\textbf{\emph{A}}_i\|^2$ and the {\em average
multiplicity} of $\lambda_j$ is
$\overline{m}_j=\frac{m_j}{n}=\|\textbf{\emph{E}}_j\|^2$. According to
(\ref{product}), this scalar product of matrices satisfies $\langle
f(\textbf{\emph{A}}),g(\textbf{\emph{A}})\rangle=\langle
f,g\rangle_\vartriangle$.

\section{Two dual approaches to almost distance-regularity}

Here we limit ourselves to the case of graphs with spectrally maximum diameter
(or the `non-degenerate' case) $D=d$. Consequently, we will use
indiscriminately the two symbols, $D$ and $d$, depending on what we are
referring to. In this context, let us consider the following two definitions of
almost distance-regularity:
\begin{defi}
\label{D1}
For a given $i$, $0\le i\le D$, a graph $\G$ is {\em $i$-punctually distance-regular} when there exist constants $p_{ji}$ such that
\begin{equation}
\label{def(a)}
\textbf{A}_i\textbf{E}_j = p_{ji}\textbf{E}_j
\end{equation}
for every $j=0,1,\ldots,d$;
 and $\G$ is {\em $m$-partially distance-regular} when it is $i$-punctually distance-regular for all $i\le m$.
 \end{defi}
 \begin{defi}
\label{D2}
For a given $j$, $0\le j\le d$, a graph $\G$ is {\em $j$-punctually eigenspace distance-regular} when
there exist constants $q_{ij}$ such that
\begin{equation}\label{def(b)}
  \textbf{E}_j\circ \textbf{A}_i = q_{ij}\textbf{A}_i
\end{equation}
for every $i=0,1,\ldots,D$;
 and $\G$ is {\em $m$-partially eigenspace distance-regular} when it is $j$-punctually eigenspace distance-regular for all $j\le m$.
\end{defi}
Notice that the concepts of $D$-partial distance-regularity and $d$-partial
\linebreak eigenspace distance-regularity coincide with the known dual
definitions of distance-regularity (see \cite{bcn}).

Some basic characterizations of punctual distance-regularity, in terms of the
distance matrices and the idempotents, were given in \cite{ddfgg10}.

\begin{pro}[\cite{ddfgg10}]
\label{first-charac(D1)} Let $D=d$. Then, $\G$ is $i$-punctually
distance-regular if and only if any of the following conditions holds:
\begin{enumerate}
 \item[$(a1)$]  $\textbf{A}_i\in {\cal A}$,
 \item[$(a2)$] $p_i(\textbf{A})\in {\cal D}$,
 \item[$(a3)$] $\textbf{A}_i=p_i(\textbf{A})$.
 \end{enumerate}
\end{pro}
Following the
duality between Definitions \ref{D1} and \ref{D2}, it seems natural to conjecture
the dual of this proposition:
A graph $\G$ is {\em $j$-punctually eigenspace
distance-regular}
if and only if any of the following conditions is satisfied:
\begin{enumerate}
\item[$(b1)$]  $\textbf{\emph{E}}_j\in {\cal D}$,
\item[$(b2)$] $q_j[\textbf{\emph{A}}]\in {\cal A}$,
\item[$(b3)$] $\textbf{\emph{E}}_j=q_j[\textbf{\emph{A}}]$,
\end{enumerate}
where $f[\textbf{\emph{A}}]=\frac{1}{n}\sum_{i=0}^d
f(\lambda_i) \textbf{\emph{A}}_i$.
However, although $(b1)$ is clearly equivalent to Definition
\ref{D2} and $(b3)\Rightarrow (b1),(b2)$, until now we have not
been able to prove any of the other equivalences and we leave
them as conjectures.

In order to derive some new characterizations of punctual distance-regularity,
besides the already defined $\overline \delta_i$ and $\overline m_j$, we
consider the following average numbers:
\begin{itemize}
\item
The {\em average crossed
local multiplicities} are
\begin{equation}\label{averagecrossed}
\overline{m}_{ij}=\frac1{n\overline{\delta}_i}\sum_{\dist(u,v)=i}m_{uv}(\lambda_j)
=\frac{\langle
\textbf{\emph{E}}_j,\textbf{\emph{A}}_i\rangle}{\|\textbf{\emph{A}}_i\|^2},
\end{equation}
where $m_{uv}(\lambda_j)=(\textbf{\emph{E}}_j)_{uv}$ are the {\em crossed local
multiplicities}.
\item
The {\em average number of shortest $i$-paths from a vertex} is
\begin{equation}
\label{meanPi}
\overline P_i= \frac{1}{n}\sum_{u\in V} P_i(u)=\frac{1}{n}\textrm{sum}\,(\textbf{\emph{A}}^i\circ \textbf{\emph{A}}_i)=\langle \textbf{\emph{A}}^i, \textbf{\emph{A}}_i\rangle=
\frac{1}{\omega_i}\langle p_i(\textbf{\emph{A}}), \textbf{\emph{A}}_i\rangle,
\end{equation}
where $P_i(u)$ denotes the number of shortest paths from a vertex $u$ to
the vertices in $\Gamma_i(u)$ and $\omega_i=(\gamma_1 \gamma_2 \cdots
\gamma_i)^{-1}$ is the leading coefficient of $p_i$, $i=1,\ldots,d$.

\item
The {\em average
number of shortest $i$-paths} is
\begin{equation}
\label{mean-aii}
\overline a_i^{(i)}=\frac{1}{n\overline\delta_i}\textrm{sum}\,(\textbf{\emph{A}}^i \circ \textbf{\emph{A}}_i)=\frac{\overline P_i}{\overline\delta_i}.
\end{equation}
\end{itemize}

\begin{pro}
\label{propo-punt-dr} Let $\G$ be a graph with predistance polynomials $p_i$ and
recurrence coefficients $\gamma_i,\alpha_i,\beta_i$,
$i=0,1,\ldots, d$. Then, $\G$ is $i$-punctually distance-regular if and only if any of the following equalities holds:
\begin{itemize}
\item[$(a1)$]
$\displaystyle \frac{1}{\overline{\delta}_i} = \sum_{j=0}^d\frac{\overline{m}_{ij}^2}{\overline m_j}$.
\item[$(a2)$]
$
\overline P_i= \frac{1}{\omega_i}\sqrt{p_i(\lambda_0)\overline \delta_i}
=\sqrt{\beta_0\beta_1\cdots \beta_{i-1}\overline \delta_i \gamma_i\gamma_{i-1}\cdots \gamma_1}$.
\item[$(a3)$]
$
\label{bound-aii}
\omega_i\overline a_i^{(i)}=1 \quad\mbox{and}\quad \overline \delta_i=p_i(\lambda_0)$.
\end{itemize}
Moreover, $\G$ is $j$-punctually eigenspace distance-regular if and only if
\begin{itemize}
\item[$(b1)$]
$\displaystyle \overline m_j=\sum_{i=0}^D\overline{\delta}_i\overline{m}_{ij}^2$.
\end{itemize}
\end{pro}
\begin{pf}
$(a1)$ This is a result from \cite{ddfgg10}.

 $(a2)$ From (\ref{meanPi}) and the Cauchy-Schwarz inequality, we get
\begin{equation}
\label{boundPi}
\omega_i \overline P_i =\langle p_i(\textbf{\emph{A}}), \textbf{\emph{A}}_i\rangle
                 \le  \|p_i(\textbf{\emph{A}})\|\|\textbf{\emph{A}}_i\|=\sqrt{p_i(\lambda_0)\overline \delta_i} \nonumber\\
                 = \sqrt{\frac{\beta_0\beta_1\cdots \beta_{i-1}}{\gamma_1\gamma_2\cdots \gamma_i}\overline \delta_i}.
\end{equation}
Moreover, equality occurs if and only if the
matrices $p_i(\textbf{\emph{A}})$ and $\textbf{\emph{A}}_i$ are
proportional, which is equivalent to $\G$ being $i$-punctually
distance-regular by Proposition \ref{first-charac(D1)}.

 $(a3)$ From (\ref{mean-aii}) and (\ref{boundPi}) we have that
 $\omega_i \overline a_i^{(i)}\le \sqrt{p_i(\lambda_0)/\overline\delta_i}$,
with equality if and only if $\G$ is $i$-punctually distance-regular. Thus, if
the conditions in $(a3)$ hold, $\G$ satisfies the claimed property. Conversely,
if $\G$ is $i$-punctually distance-regular, both equalities in $(a3)$ are
simple consequences of $p_i(\textbf{\emph{A}})=\textbf{\emph{A}}_i$. Indeed,
the first one comes from considering the $uv$-entries, with $\dist(u,v)=i$, in
the above matrix equation, whereas the second one is obtained by taking square
norms.

$(b1)$ From (\ref{averagecrossed}), we find that the orthogonal projection of
$\textbf{\emph{E}}_j$ on ${\cal D}$ is $ \widehat{\textbf{\emph{E}}_j}
=\sum_{i=0}^D \overline m_{ij}\textbf{\emph{A}}_i $. Now, from
$\|\widehat{\textbf{\emph{E}}_j}\|^2\le \|\textbf{\emph{E}}_j\|^2$ we get
 $$
\sum_{i=0}^D \overline m_{ij}^2\|\textbf{\emph{A}}_i\|^2
 = \sum_{i=0}^D\overline{\delta}_i\overline{m}_{ij}^2\le \overline{m}_j
 $$
and, in the case of equality, Definition \ref{D2} applies with
$q_{ij}=\overline m_{ij}$.
\end{pf}

Notice the duality between $(a1)$ and $(b1)$ with $\frac{1}{\overline{\delta}_i}$ and $\overline{m}_j$.

Now, let us consider the more global concept of partial distance-regularity. In
this case, we also have the following new result where, for a given $0\le i\le
d$, $s_i=\sum_{j=0}^i p_j$, $t_i=H-s_{i-1}=\sum_{j=i}^d p_j$,
$\textbf{\emph{S}}_i=\sum_{j=0}^i \textbf{\emph{A}}_j$, and
$\textbf{\emph{T}}_i=\textbf{\emph{J}}-\textbf{\emph{S}}_{i-1}=\sum_{j=i}^d
\textbf{\emph{A}}_j$.
\begin{pro}
\label{charac-pardr}
A graph $\G$ is $m$-partially distance-regular if and only if any of the following conditions holds:
\begin{itemize}
\item[$(a1)$]
$\G$ is $i$-punctually distance-regular for $i=m,m-1,\ldots,\emph{max}\{2,2m-d\}$.
\item[$(a2)$] $\G$ is $m$-punctually distance-regular and
    $t_{m+1}(\textbf{A})\circ \textbf{S}_m=\textbf{O}$.
\item[$(a3)$]
$s_i(\textbf{A})=\textbf{S}_i$ for $i=m,m-1$.
\end{itemize}
\end{pro}
\begin{pf}
In all cases, the necessity is clear since
$p_i(\textbf{\emph{A}})=\textbf{\emph{A}}_i$ for every $0\le i\le m$ (for
$(a2)$, note that
$t_{m+1}(\textbf{\emph{A}})=\textbf{\emph{J}}-s_m(\textbf{\emph{A}})$). Then,
let us prove sufficiency. The result in $(a1)$ is basically Proposition 3.7 in
\cite{ddfgg10}. In order to prove $(a2)$, we show by (backward) induction that
$p_i(\textbf{\emph{A}})=\textbf{\emph{A}}_i$ and
$t_{i+1}(\textbf{\emph{A}})\circ \textbf{\emph{S}}_i=\textbf{\emph{O}}$ for $i=
m,m-1,...,0.$ By assumption, these equations are valid for $i=m$. Suppose now
that $p_i(\textbf{\emph{A}})=\textbf{\emph{A}}_i$ and
$t_{i+1}(\textbf{\emph{A}})\circ \textbf{\emph{S}}_i=\textbf{\emph{O}}$ for
some $i>0$. Then, $t_i(\textbf{\emph{A}})\circ
\textbf{\emph{S}}_i=\textbf{\emph{A}}_i$ and, multiplying both terms by
$\textbf{\emph{S}}_{i-1}$ (with the Hadamard product), we get
$t_{i}(\textbf{\emph{A}})\circ \textbf{\emph{S}}_{i-1}=\textbf{\emph{O}}$. So,
what remains is to show that
$p_{i-1}(\textbf{\emph{A}})=\textbf{\emph{A}}_{i-1}$. To this end, let us
consider the following three cases:
\begin{itemize}
\item[$(i)$] For $\dist(u,v)>i-1$, we have
    $(p_{i-1}(\textbf{\emph{A}}))_{uv}=0$.
\item[$(ii)$] For $\dist(u,v)=i-1$, we have
    $(t_{i+1}(\textbf{\emph{A}}))_{uv}=0$, so
    $(p_{i-1}(\textbf{\emph{A}}))_{uv}=(s_{i-1}(\textbf{\emph{A}}))_{uv}$ $=
    (s_{i-1}(\textbf{\emph{A}}))_{uv}+(\textbf{\emph{A}}_{i})_{uv}=(s_{i}(\textbf{\emph{A}}))_{uv}=
  1-(t_{i+1}(\textbf{\emph{A}}))_{uv}=1$.
\item[$(iii)$] For $\dist(u,v)<i-1$, we use the recurrence
    (\ref{recur-pol}) to write
\begin{eqnarray*}
  xt_i=\sum_{j=i}^d xp_j & = & \sum_{j=i}^d (\beta_{j-1}p_{j-1}+\alpha_j p_j+ \gamma_{j+1}p_{j+1})\\
  & = & \beta_{i-1}p_{i-1}- \gamma_i p_i + \sum_{j=i}^d(\alpha_j+\beta_j+\gamma_j) p_j \\
   &  = & \beta_{i-1}p_{i-1}- \gamma_ip_i+ \delta t_i ,
\end{eqnarray*}
which gives
$$
\textbf{\emph{A}}t_i(\textbf{\emph{A}})=\beta_{i-1}p_{i-1}(\textbf{\emph{A}}) -\gamma_i \textbf{\emph{A}}_i+ \delta t_i(\textbf{\emph{A}}).
$$
Then, since
$(t_i(\textbf{\emph{A}}))_{uv}=(\textbf{\emph{A}}_i)_{uv}=0$
and $\beta_{i-1}\neq 0$, we get
$$
(p_{i-1}(\textbf{\emph{A}}))_{uv}=\frac1{\beta_{i-1}}(\textbf{\emph{A}}t_i(\textbf{\emph{A}}))_{uv}=
\frac1{\beta_{i-1}}\sum_{w\in
\Gamma(u)}(t_i(\textbf{\emph{A}}))_{wv}=0,
$$
because $\dist(v,w)\le \dist(v,u)+\dist(u,w) \le i-1$ for
the relevant $w$.
\end{itemize}
From  $(i),(ii),$ and $(iii)$, we have that
$p_{i-1}(\textbf{\emph{A}})=\textbf{\emph{A}}_{i-1}$, so by
induction $\G$ is $m$-partially distance-regular, and the
sufficiency of $(a2)$ is proven. Finally, the sufficiency of $(a3)$
follows from that of $(a2)$ because
$s_i(\textbf{\emph{A}})=\textbf{\emph{S}}_i$ for every $i\in\{m-1,m\}$
implies that
$p_m(\textbf{\emph{A}})=(s_m-s_{m-1})(\textbf{\emph{A}})=\textbf{\emph{S}}_m-\textbf{\emph{S}}_{m-1}=\textbf{\emph{A}}_m$
and $t_{m+1}(\textbf{\emph{A}})\circ
\textbf{\emph{S}}_m=(\textbf{\emph{J}}-s_m(\textbf{\emph{A}}))\circ
\textbf{\emph{S}}_m=(\textbf{\emph{J}}-\textbf{\emph{S}}_m)\circ
\textbf{\emph{S}}_m= \textbf{\emph{O}}$.
\end{pf}

Given some vertex $u$ and an integer  $i\le \textrm{ecc}(u)$, we denote by
$N_i(u)$ the {\em $i$-neighborhood} of $u$, which is the set of vertices that
are at distance at most $i$ from $u$. In \cite{f02} it was proved that
$s_i(\lambda_0)$ is upper bounded by the harmonic mean of the numbers
$|N_i(u)|$ and equality is attained if and only if
$s_i(\textbf{\emph{A}})=\textbf{\emph{S}}_i$. A direct consequence of this
property and Proposition \ref{charac-pardr}$(a3)$ is the following
characterization.

\begin{thm}
\label{thm-pdr}
A graph $\G$ is $m$-partially distance-regular if and only if, for every $i\in \{m-1,m\}$,
$$
s_i(\lambda_0) = \frac{n}{\sum_{u\in V}|N_i(u)|^{-1}}.
$$
\end{thm}

\section{Distance-regular graphs}
Let us particularize our results to the case of distance-regular graphs. With
this aim, we use the following theorem giving some known  characterizations.

\begin{thm}[\cite{f01,fgy1b}]
\label{charac-drg}  A graph
$\G$  with $d+1$ distinct eigenvalues and diameter $D=d$ is
distance-regular if and only if any of the following statements
is satisfied:
\begin{itemize}
\item[$(a)$] $\G$ is $D$-punctually distance-regular.
\item[$(b)$] $\G$ is $j$-punctually eigenspace distance-regular for
    $j=1,d$.
\end{itemize}
\end{thm}
In fact, notice that $(a)$ corresponds to any of the conditions in Proposition \ref{charac-pardr} with $m=d$.
Moreover, the duality between $(a)$ and $(b)$ is made apparent when they are
stated as follows:
\begin{itemize}
\item[$(a)$] $\textbf{\emph{A}}_0(=\textbf{\emph{I}}),\textbf{\emph{A}}_1(=\textbf{\emph{A}}),\textbf{\emph{A}}_D\in {\cal A}$;
\item[$(b)$] $\textbf{\emph{E}}_0(=\frac{1}{n}\textbf{\emph{J}}),\textbf{\emph{E}}_1,\textbf{\emph{E}}_d\in {\cal D}$.
\end{itemize}

Then, by using Theorem \ref{charac-drg} and  Proposition
\ref{propo-punt-dr}$(a1)$ and $(b1)$, and Theorem \ref{thm-pdr} (with $m=d$),
we have the spectral excess theorem \cite{fg97} in the next condition $(a)$,
its dual form in $(b)$, and its harmonic mean version \cite{f02,vd08} in $(c)$.

\begin{thm}
A regular graph $\G$ with $D=d$ is distance-regular if and only if any of the
following equalities holds:
\begin{itemize}
\item[$(a)$]
$\displaystyle \frac{1}{\overline{\delta}_d}= \sum_{j=0}^d\frac{\overline{m}_{dj}^2}{\overline m_j}$.
\item[$(b)$]
$\displaystyle  \overline m_j=\sum_{i=0}^D\overline{\delta}_i\overline{m}_{ij}^2\ \mbox{ for $j=1,d$}$.
\item[$(c)$] $\displaystyle s_{d-1}(\lambda_0) = \frac{n}{\sum_{u\in V}|N_{d-1}(u)|^{-1}}$.
\end{itemize}
\end{thm}

In fact, condition $(a)$ is usually written in its equivalent form
$\overline{\delta}_d=p_d(\lambda_0)$ as, when $i=d$, the first condition in
Proposition \ref{bound-aii}$(a.3)$ always holds since
$$
\overline a_d^{(d)}=\frac{1}{\overline \delta_d}\langle\textbf{\emph{A}}^d, \textbf{\emph{A}}_d\rangle=
\frac{1}{\overline \delta_d \omega_d}\langle H(\textbf{\emph{A}}), \textbf{\emph{A}}_d\rangle=
\frac{1}{\overline \delta_d \omega_d}\langle \textbf{\emph{J}}, \textbf{\emph{A}}_d\rangle=
\frac{1}{\overline \delta_d \omega_d}\| \textbf{\emph{A}}_d\|^2=\frac{1}{\omega_d}.
$$
Notice also that, in $(c)$, we do not need to impose the condition of Theorem \ref{thm-pdr} for $i=d$  since $s_d(\lambda_0)=H(\lambda_0)=N_d(u)=n$ for every $u\in V$.


\begin{thebibliography}{00}
{\small
\bibitem{biggs}
N. Biggs, {\it Algebraic Graph Theory}, Cambridge University Press,
Cambridge, 1974, second edition, 1993.

\bibitem{bcn}
A.E. Brouwer, A.M. Cohen, and A. Neumaier, {\it Distance-Regular
Graphs}, Springer-Verlag, Berlin-New York, 1989.

\bibitem{cffg09}
M. C\'amara, J. F\`abrega, M.A. Fiol, and E. Garriga,
Some families of orthogonal polynomials of a discrete variable and
their applications to graphs and codes, {\em Electron. J. Combin.} {\bf 16(1)} (2009), \#R83.

\bibitem{cds80} C.D. Cvetkovi\'c, M. Doob, H. Sachs, {\it Spectra of Graphs},
    third edition, Johann Barth Verlag, 1995. First edition: Deutscher Verlag
    der Wissenschaften, Academic Press, Berlin, New York, 1980.

\bibitem{ddfgg10} C. Dalf\'{o}, E.R. van Dam, M.A. Fiol, E. Garriga, and B.L.
    Gorissen, On almost distance-regular graphs, {\em J. Combin. Theory Ser.
    A} {\bf 118} (2011), 1094--1113.

\bibitem{dfg09}
C. Dalf\'{o}, M.A. Fiol, and E. Garriga, On
$k$-walk-regular graphs, \emph{Electron. J. Combin.} {\bf 16(1)} (2009), \#R47.

\bibitem{f01} M.A. Fiol, On pseudo-distance-regularity, {\it Linear Algebra
    Appl.} {\bf 323} (2001), 145--165.

\bibitem{f02}
M.A. Fiol, Algebraic characterizations of
distance-regular graphs, {\it Discrete Math.} {\bf 246}
(2002), 111--129.

\bibitem{fg97} M.A. Fiol and E. Garriga, From local adjacency
polynomials to locally pseudo-distance-regular graphs, {\it
J.  Combin. Theory Ser. B} {\bf 71} (1997), 162--183.

\bibitem{fgg09}
M.A. Fiol, S. Gago, and E. Garriga, A simple proof of the spectral excess theorem for
distance-regular graphs, {\it Linear Algebra Appl.} {\bf 432} (2010), 2418--2422.

\bibitem{fgy1b} M.A. Fiol, E. Garriga, and J.L.A. Yebra, Locally
    pseudo-distance-regular graphs, {\it J.  Combin. Theory Ser. B} {\bf 68}
    (1996), 179--205.

\bibitem{g93}
C.D. Godsil, {\it Algebraic Combinatorics}, Chapman and Hall, NewYork, 1993.

\bibitem{hof63}
A.J. Hoffman, On the polynomial of a graph,
{\it Amer. Math. Monthly} {\bf 70} (1963), 30--36.

\bibitem{p91} D.L. Powers,  Partially distance-regular graphs,  in  {\em Graph
    Theory, Combinatorics, and Applications, Vol. 2. Proc. Sixth Quadrennial
    Int. Conf. on the Theory and Appl. of Graphs, Western Michigan University,
    Kalamazoo, 1988} (Y. Alavi et al., eds.), Wiley, New York, 1991, 991--1000.

\bibitem{vd08} E.R. van Dam, The spectral excess theorem for distance-regular
    graphs: a global (over)view, {\em Electron. J. Combin.} {\bf 15(1)} (2008),
    \#R129.}
\end{thebibliography}
\end{document}